\documentclass[12pt]{article}
\usepackage{latexsym,amssymb,amscd,amstext}
\usepackage[all]{xy}
\begin{document}
\title{Relativity groupoid, instead of relativity group}
\author{Zbigniew Oziewicz\\
Universidad Nacional Aut\'onoma de M\'exico\\Facultad de Estudios Superiores Cuautitl\'an\\
Apartado Postal \# 25, C.P. 54714 Cuautitl\'an Izcalli\\Estado de M\'exico, M\'exico\\
oziewicz@servidor.unam.mx}\date{Submitted 20 June 2006, Accepted 18
July 2006\\Revised 28 April 2007}\maketitle

\newtheorem{groupoid}{Definition}[section]
\newtheorem{odef}[groupoid]{Definition}
\newtheorem{oaxiom}[groupoid]{Axiom}
\newtheorem{assum}[groupoid]{Assumptions}
\newtheorem{ocor}[groupoid]{Corollary}
\newtheorem{othm}[groupoid]{Theorem}
\newtheorem{olem}[groupoid]{Lemma}

\newcommand{\C}{\mathbb{C}}\newcommand{\Real}{\mathbb{R}}
\newcommand{\N}{\mathbb{N}}\newcommand{\Qr}{\mathbb{Q}}
\newcommand{\La}{\mathcal{L}}\newcommand{\F}{\mathcal{F}}
\newcommand{\ie}{\textit{i.e.}\,}
\newcommand{\alg}{\text{alg}}
\newcommand{\obj}{\text{obj}}
\newcommand{\arrows}{\text{arrows}}
\newcommand{\cat}{\text{cat}}
\newcommand{\trace}{\text{tr}}
\newcommand{\im}{\text{im}}
\newcommand{\End}{\text{End\,}}
\newcommand{\Span}{\text{span}}
\newcommand{\der}{\text{der}}
\newcommand{\obs}{\text{Obs}}
\newcommand{\vel}{\mathbf{v}}\newcommand{\cel}{\mathbf{c}}
\newcommand{\uel}{\mathbf{u}}\newcommand{\xel}{\mathbf{x}}\newcommand{\yel}{\mathbf{y}}
\newcommand{\wel}{\mathbf{w}}\newcommand{\pel}{\mathbf{p}}
\newcommand{\zel}{\mathbf{0}}\newcommand{\opp}{\text{opp}}
\newcommand{\vvel}{\textstyle{\frac{\vel}{c}}}

\begin{abstract} Minkowski in 1908 used space-like \textit{binary} velocity-field of a medium, relative to an observer.
Hestenes in 1974 introduced, within a Clifford algebra, an axiomatic binary relative velocity as a Minkowski bivector.
We propose consider binary relative velocity as a traceless nilpotent endomorphism in an operator algebra.
Any concept of a \textit{binary} axiomatic relative velocity made possible the replacement of the Lorentz
relativity group by the relativity groupoid. The relativity groupoid is a category of massive bodies in mutual relative motions,
where a binary relative velocity is interpreted as a categorical morphism with the associative addition.
This associative addition is to be contrasted with non-associative addition of (ternary) relative velocities in
isometric special relativity (loop structure).

We consider an algebra of many time-plus-space splits, as an operator algebra generated by idempotents.
The kinematics of relativity groupoid is ruled by associative Frobenius operator algebra,
whereas the dynamics of categorical relativity needs the non-associative Fr\"olicher-Richardson operator algebra.

The Lorentz covariance is the cornerstone of physical theory.
Observer-dependence within relativity groupoid, and the
Lo\-rentz\--co\-var\-iance (within the Lorentz relativity group),
are different concepts. Laws of Physics could be observer-free,
rather than Lorentz-invariant.\end{abstract}

\noindent\textbf{Keywords.} Groupoid category; binary
velocity-morphism; associative addition of binary relative
velocities.

\markboth{Zbigniew Oziewicz}{Relativity Groupoid instead of
Relativity Group}

\section{Galileo versus Descartes}

Aether was introduced by Ren\'e Descartes around 1650. Descartes conceive the space as an \textit{absolute substance,} substratum, medium, that exists independently of any matter in it. Like a hall with sitting chairs-locations waiting for spectators. Galileo claim instead, in 1632, that space is \textit{relative}. The relativity
theory is about the relativity of proper-time-plus-space splits, and not about the change of numerical coordinates.

Why does a space, a set of locations, not need be neither the absolute nor primitive concept, as it was postulated by Descartes, and also by Newton in 1686? How could it be \textit{relative}, \ie dependent on another primitive concept? What it is this primitive concept that space and time are derived? Galileo conceived space to be massive bodies-dependent. Galilean space has no reality without the bodies
that `it contains'. Galilean primitive concept from which space is
derived is the relative \textit{motion}, the relative \textit{velocity} among massive bodies. Because of relative velocities, there are a multitude of relative spaces. If all massive bodies are at rest relative to each other, zero
relative velocities only, or just one massive body, then space would be
absolute and proper-time (a set of instants) would also be absolute.
Why is space relative? Because there are relative velocities, and
hence the systems of more than just one massive body in mutual motions (Galileo 1632).

In this way we come to a vicious circle: in order to define a velocity it is said, according to Newton, that we need first a primitive space and a primitive time, and a velocity is a derivative of a space with respect to a time. However, the very concept of a space (and a time), according to Galileo, is very much relative
because of multitude of a priori relative velocities. Therefore one can ask what comes first: chicken or egg? space \& time or relative velocities among massive bodies? \'Emil Picard conceived motion and time to be dual: \textit{time is measuring a motion, and a motion is measuring a time}.

Every velocity is relative, it is a velocity of one massive body relative to another massive body. Within the isometric special relativity the relative velocity is defined as a parameter in the Lorentz boost. This entails relative velocity to be \textit{ternary}. A ternary relative velocity is a velocity among two massive bodies, as measured/seen by the third preferred exterior observer (Oziewicz 2005). Minkowski in 1908 used \textit{binary} velocity field of a medium relative to observer. Hestenes in 1974, within Clifford geometric-algebra, introduced basis-free and coordinate-free binary relative velocity as a Minkowski bivector. Let an observer be given by a Minkowski time-like vector field $P.$ Hestenes defined the physical space of an observer $P$ as a subspace of all Minkowski simple bivectors $\{b\},$ such that vanishes the Grassmann tri-vector, $P\wedge b=0\;\Leftrightarrow\;b=(i_{gP}b)\wedge P.$ For the Hestenes formalism we refer to (Hestenes 2003; Baylis with Sobczyk 2004). We define the physical space of an observer $P$ to be orthogonal to a vector field $P,$ \ie a sub-module of all Minkowski's space-like vector-fields that are orthogonal to $P,$ $P\cdot\vel=0\;\Leftrightarrow\;\vel=i_{gP}b$ (here a bivector $b$ needs not to be simple). For a simple bivectors this definition is in a bijective correspondence with the Hestenes physical space of Minkowski's simple bivectors.

The set of all relative velocities among massive bodies is \textit{not} a vector space of linear algebra because not every pair of such relative velocities can be composed (velocity addition is the partial operation), nor is the commutativity of composition (an additive Abelian group structure) applicable when it is defined.
This holds equally well in Newtonian physics with absolute simultaneity (identified with the physical time), and as well as for relative simultaneity and finite light velocity. Moreover, every massive body possess its own identity velocity-morphism, it is the zero velocity of this object relative to himself. There is no universal unique zero velocity that would be massive-body-independent. The zero velocity of the Earth relative to Earth must \textit{not} be identified with the zero velocity of the Sun relative to the Sun.

If relative velocity is not a vector of linear algebra, then one can ask with which mathematical concept the physical relative velocity could be identified? Each binary relative velocity is a binary function of ordered pair of massive bodies, like a set-valued Hom functor in category theory.

\newpage
\section{Groupoid category}

A category consists of a family of objects and a family of arrows/morphisms. Every arrow has a source object and a target object. Two arrows are said to be composable if the target object of one arrow is the same as the source object of the second arrow. The composition of composable arrows is assumed to be associative.

\begin{groupoid}[Groupoid category] \emph{A category is said to be a \textit{groupoid category}, if and only if every morphism has a two-sided inverse. In particular a \textit{group} is a groupoid one-object-category, with just one object, hence with universal unique neutral element-morphism. A groupoid category is said to be \textit{transitive, or connected,} if there is an arrow joining any two of its objects.}\end{groupoid}

\begin{odef}[Terminal and initial object] \emph{An object $q$ is terminal if to each object $p$ there is exactly one arrow $p\rightarrow q.$ An object $p$ is initial if to each object $q$ there is exactly one arrow $p\rightarrow q.$ A \textit{null} object is an object which is both initial and terminal.}\end{odef}

If $p$ is terminal, the only arrow $p\rightarrow p$ is the identity. Any two terminal objects are isomorphic (Mac Lane 1998, p. 20, p. 194). A transitive groupoid category where \textit{every} object is null, is called a \textit{pair groupoid}.

\section{Relativity Groupoid as a Pair Groupoid}

Let a base of a groupoid be given by a set of massive bodies (not necessarily inertial), and let each morphism be the unique binary relative velocity. Such pair groupoid of null object is said to be the relativity groupoid,
\begin{equation}\varpi=\{\obj\,\varpi,\arrows\,\varpi\}=\{\text{massive bodies},
\text{binary relative velocities}\}.\end{equation}

Unique arrow from an object $p$ to an object $q$ is denoted by $\varpi(p,q),$ with analogy to Hom-set bifunctor. For each object $p\in\obj\,\varpi,$ a map
\begin{equation}\begin{CD}\obj\,\varpi\ni q\quad@>{\varpi_p}>>\quad\varpi(p,q)
\in\arrows\,\varpi,\end{CD}\end{equation} extends to covariant representable functor among connected groupoid categories. The object $p$ is representing-object for a functor $\varpi_p.$ We call this structure the connected groupoid (1,1)-category, or enriched pair-groupoid category. The relativity groupoid is (arrows\,$\varpi$)-enriched groupoid category, rather then more restrictive concept of 2-category. This category is neither Abelian, nor additive.

Each relative binary velocity is a categorical morphism (arrow that need not be a map because objects need not to be sets having elements). An arrow $\varpi(p,q)$ is interpreted as a binary velocity-morphism of a massive body $q\in\obj\varpi$ relative to a body $p\in\obj\varpi,$ \ie a velocity as measured by $p.$ We say that the source (or domain) of a velocity-morphism $\varpi(p,q)$ is $p,$ and the target (or codomain) of $\varpi(p,q)$ is a body $q.$ A category symbol $\varpi,$ is
interpreted as an actual binary-velocity-measuring device. We display $\varpi(p,q)$ as a categorical arrow (morphism, directed-path) which originates (is outgoing) at observer $p,$ $p$ is a node of the directed graph, and terminates (is incoming) at an observed body $q,$
\begin{equation}\xymatrix{\ldots\ar[r]&{\quad p\quad}\ar@/^1pc/[rr]^{\varpi(p,q)}&&
{\quad q}\ar@/^1pc/[ll]_{\varpi(q,p)}\quad\ar[r]&\ldots;}\label{display}
\end{equation}\begin{tabular}{cc}Observer of $\varpi(p,q)$ is $p.$&
Observed body with $\varpi(p,q)$ is $q.$\\
Observed body with $\varpi(q,p)$ is $p.$&Observer of $\varpi(q,p)$ is $q.$
\end{tabular}\medskip

There is not yet, either the concept of a spacetime, or the concept of a time, or space. The relative binary-velocity-morphism, and massive indivisible objects, are primary, postulated primitive concepts. Any massive body (observer, ob\-ser\-ved), is a null object in a category $\varpi.$ The categorical null object is indivisible, like the Leibniz monad, it is not a set. The binary relative-velocity is the \textit{primitive} notion, and we are introducing this notion as the \textit{morphism} in the pair-groupoid category of abstract observers. Each morphism of the pair-groupoid category $\varpi$ is a binary (interior) relative velocity; the unique binary velocity-morphism instead of a set of the isometric Lorentz transformations in special relativity (with non-unique, ternary-velocities). We are not associating the concept of the binary relative-velocity either with the Lorentz boost, or with the concept of the vector space. We wish to see the relative-velocity exclusively as the categorical morphism in pair-groupoid category that is not Abelian. What is relativity groupoid? It is a family of all binary relative-velocities.
\begin{figure}[h]{\Large$$\xymatrix{&q\ar@/^4pc/[dr]^{\varpi(q,r)}
\ar@/^/[dl]_{\varpi(q,p)}\\
p\ar@/^4pc/[ur]^{\varpi(p,q)}\ar@/_2pc/[rr]_{\varpi(p,r)}&&
r\ar@/^/[ul]_{\varpi(r,q)}\ar@/^/[ll]_{\varpi(r,p)}}$$}\caption{A relativity groupoid $\varpi$ of three massive bodies, $p,q,r,$ with six binary velocities-morphisms. This category generate 9-dimensional algebra. The identity arrows are suppressed.\label{category}}\end{figure}
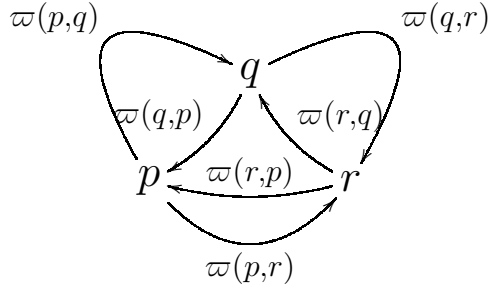
In relativity groupoid there is no need to distinguish separately the constant velocities (special relativity) from the variable accelerated velocities, hence the relativity groupoid goes beyond boundary of the special relativity.

\section{Relativity Groupoid is an Algebra}\setcounter{equation}{3}

Let $\F$ be an associative, unital and commutative ring. We denote by $\Span_\F\varpi$ the $\F$-module with a category $\varpi$ as a set of basic vectors. This $\F$-module $\Span_\F\varpi$ consists of all formal $\F$-linear combinations of the elements of $\varpi,$ \ie the formal combinations that mix objects with arrows.

It is postulated that an $\F$-module $\Span_\F\varpi$ is an associative $\F$-algebra, called an algebra `of massive bodies/observers in the mutual relative motions', and denoted by $\obs(\varpi).$ The algebra $\obs(\varpi)$ is generated by
(presented on) objects and arrows of $\varpi,$ subject to the relations: It is postulated that every object $p\in\obj\varpi$ is an idempotent, $p^2=p\in\obs(\varpi).$ Every object of $\varpi,$ seen in an algebra $\obs(\varpi),$ and represented as an operator in $\End(\der\F),$ looks like a pure state in quantum mechanics.

The above postulate have the following motivation. In order to have just one space, we need to fix one massive body as the reference system. The correspondence, massive body $\leftrightarrow$ idempotent, is motivated by a desire that every massive body $p\in\End(\der\F),$ splits \begin{equation}\der\F=(\ker p)\oplus(\im p)
=(\text{space})\oplus(\text{time}).\end{equation} The choice of one body, for example the Earth, as the reference system, does \textit{not} require coordinates. Such choice is coordinate-free, and basis-free. We call any massive body, the Earth, the Moon, an observer (no measuring devices, rods and clocks are involved).

An $\F$-dual $\F$-module is denoted by $(\Span_\F\varpi)^*.$ There is the distinguished covector $\trace\in(\Span_\F\varpi)^*,$ $\trace:\Span_\F\varpi\rightarrow\F.$ If $p$ is an object of $\varpi,$ then $\trace\,p=1.$

An arrow $\varpi(p,q),$ outgoing at $p,$ and incoming at $q,$ is given by nilpotent traceless endomorphism,
\begin{equation}\frac{\varpi(p,q)}{c}\;\equiv\;\frac{qp}{\trace(qp)}-p\quad\Longrightarrow\quad (\varpi(p,q))^2=0\in\obs(\varpi),\quad\trace\,\varpi(p,q)=0.\label{val}\end{equation}
Therefore the covector $\trace$ is distinguishing an object $p,$ from the identity arrow $\varpi(p,p).$ For every non empty word (string) of objects, $p,q,r,\ldots,s\in\obj\varpi,$ the following properties and relation are postulated,
\begin{eqnarray}1\leq\{\trace(pqr\ldots s)\}^2=\trace(pq)
\trace(qr)\trace(r..)\ldots\trace(..s)\trace(sp),\\
1\leq\trace(pq)=\trace(qp),\quad\trace(pqr)=\trace(qpr),\label{m1}\\
\forall\;A\in\obs(\varpi),\quad
qAp=\trace(qAp)\left\{p+\frac{1}{c}\varpi(p,q)\right\},\quad
c<\infty.\label{m}\end{eqnarray}

From relations (\ref{m1}-\ref{m}), one can deduce the algebra multiplication table, the multiplication of arrows $\{\varpi(p,q)\}$ with objects $\{p,q,r,s,...\},$ and arrows with arrows (we set $c=1$),
\begin{equation}qp=\trace(qp)\{p+\varpi(p,q)\},\label{t1}\end{equation}
\begin{equation}q\varpi(p,r)=\left(\frac{\trace(qrp)}{\trace(rp)}-\trace(pq)\right)
\{p+\varpi(p,q)\},\label{t2}\end{equation}\begin{equation}
\varpi(q,r)p=\frac{\trace(rqp)}{\trace(rq)}\{p+\varpi(p,r)\}-\trace(qp)
\{p+\varpi(p,q)\},\end{equation}\begin{eqnarray}
\varpi(p,q)\varpi(r,s)=\frac{1}{\trace(qp)}\left(\frac{\trace(qpsr)}
{\trace(sr)}-\trace(qpr)\right)\{r+\varpi(r,q)\}\qquad\nonumber\\
\hspace{3cm}+\left(\trace(pr)-\frac{\trace(psr)}{\trace(sr)}\right)\{r+\varpi(r,p)\}.
\end{eqnarray}
Zero velocity of a body $p$ relative to himself is $\varpi(p,p)=\zel_p.$ All these zero velocities are equal to unique zero of algebra $\obs(\varpi)$ for all bodies, however they \textit{cannot} be identified with respect to associative composition
of morphisms in a category $\varpi,$ see Sec. 7.

In Cruz \& Oziewicz (2006) we consider the augmented unital algebra,
$(\obs\,\varpi)\oplus\F,$ and pose a hypothesis that this trace-class algebra is a Frobenius algebra for each cardinality of a finite subfamily of objects.

\section{Scalar Magnitude of Arrow}\setcounter{equation}{12}

The arrows of the relativity groupoid $\varpi$ does not possess a linear structure (an arrow multiplied by a scalar is not an arrow). An Euclidean angle between arrows with the same source only, and a scalar magnitude of each arrow, needs the following separate postulate
\begin{equation}\frac{\varpi(q,p)\cdot\varpi(q,r)}{c^2}\equiv 1-
\frac{\trace(pr)}{\trace(pqr)}\qquad\Longrightarrow\qquad\left(\frac{|\varpi(p,q)|}
{c}\right)^2=1-\frac{1}{\trace(pq)}.\label{angle}\end{equation}

\section{Module over an Algebra of Observers}\setcounter{equation}{13}

An $\F$-algebra of observers, $\obs(\varpi),$ is abstractly isolated from the concepts of spacetime, time and space. The indivisible objects of $\varpi$ are not yet located neither in a space, nor in a time. A ring $\F$ is interpreted as an $\Real$-algebra of the classical measurements, in case that $\F$ is commutative.

The spacetime manifold of events is usually identified with $\alg(\F,\Real).$ We propose to see the concept of the spacetime encoded in an $\obs(\varpi)$-module, \ie in an algebra morphism, from an associative $\F$-algebra $\obs(\varpi),$ into
endomorphism algebra of the Lie $\F$-module $\der\F,$
\begin{equation}\begin{CD}\obs(\varpi)\quad
@>{\text{associative algebra}}>{\text{morphism}}>
\quad\End(\der\F)\simeq(\der\F)\otimes(\der\F)^*.\end{CD}\end{equation}
The Lie $\F$-module $\der\F$ of derivations of the ring $\F,$ and the dual $\F$-module, $(\der\F)^*,$ of the differential forms, are considered as $\obs(\varpi)$-modules.

A massive body $p\in\obj\varpi$ is represented by $(1,1^*)$-tensor field $p\in\End(\der\F),$ where $p^2=p$ must be the minimal polynomial. Every object of $\varpi,$ seen in an algebras, $\obs(\varpi)$ and $\End(\der\F),$ looks like a pure state in quantum mechanics.

A pull back (transpose) of $p$ is denoted by $p^*,$ it is a $(1^*,1)$-tensor field. The gravity potential tensor field $g=g^*\in(\der\F)^*\otimes(\der\F)^*,$ and his inverse $g^{-1}\in(\der\F)\otimes(\der\F),$ are considered as an invertible Grassmann $\F$-algebra morphisms from multi-vector fields $(\der\F)^\wedge,$ to differential multi-forms,
\begin{equation}\xymatrix{{n}\ar[rr]^{p}&&{\quad n}\ar@/^1pc/[rr]^{g}&&{n^*\quad}
\ar@/^1pc/[ll]^{g^{-1}}\ar[rr]^{p^*}&&{n^*,}}\quad\forall\,n\in\N.\end{equation}
In this convention, $1\in\N$ denotes (a grade of) a Lie $\F$-module of vector fields $\der\F,$ $2\in\N$ denotes an $\F$-module of bivector fields, $(\der\F)\wedge(\der\F),$ $1^*$ is a module of Pfaffian differential one-forms $(\der\F)^*,$ etc. Therefore $g\in\alg(n^*,n),$ however, $\obj(\varpi)\ni
p\leadsto p\in\der(n,n).$ Here $\leadsto$ means `extends', and
$\der(n,n)$ is short for $\der((\der\F)^{\wedge n},(\der\F)^{\wedge n}).$

The composition, $g\circ p,$ and their transpose, $(g\circ p)^*=p^*\circ g,$ are both morphism with the same domain and codomain, and therefore one can ask that a massive body $p\in\End(\der\F)$ is metric-compatible ($g$-orthogonal),
\begin{equation}\dim\im\,p=1\;\&\;\trace\,p=1\quad\&\quad p^*\circ g
=g\circ p\quad\Longleftrightarrow \quad p=\frac{P\otimes(gP)}{g(P\otimes P)},\quad p^2=p.\label{grav}\end{equation}

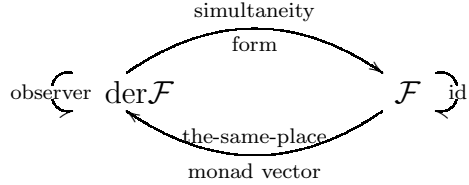
\begin{figure}[h]$$\quad\xymatrix{{}\ar@(ul,dl)[]|{\text{observer}}
\hspace{1cm}
{\quad\der\F\quad}\ar@/^2pc/[rr]^{\text{simultaneity}}_{\text{form}}&&
{\quad\F\quad}\ar@(ur,dr)[]|{\text{id}}\ar@/^2pc/[ll]^{\text{monad
vector}}_{\text{the-same-place}}&&{}}$$$$p=(\text{monad
vector})\otimes(\text{simultaneity form})\quad$$\caption{An observer
as the split-idempotent.}\end{figure}

The above compatibility of observer $p$ with the gravitational potential $g$ (\ref{grav}), it seems that could be tested experimentally rather than postulated a priori. Let $g=g^*$ be a Lorentzian metric tensor field of signature $(-+++),$ not necessarily curvature-free. One can set all monad fields (the eigenvectors of observers that belong to eigenvalue 1) to be normalized, $pP=P$ with $P^2=-1.$ For more motivations about the concept of massive observer as $(1,1)$-tensor field, we refer to Kocik (1997), and Cruz \& Oziewicz (2003). Essentially such observer-idempotent is the same as non-viscous relativistic fluid given by energy-momentum endomorphism field with mass density and with pressure (Tolman 1918, 1987). The space-like simultaneity of an observer $p$ is given by the Einstein-Minkowski proper time differential one-form, $gP$ (Minkowski 1908). This time-like differential Pfaff form, $gP,$ encode the unique
empirical and metric-dependent simultaneity of an observer $p.$

Electro-magnetics of moving bodies within relativity groupoid is different from the isometric special relativity with Lorentz transformations, however this goes beyond the scope of this note.

\section{Associative Addition of Binary Relative-Velocities}
\setcounter{equation}{16}

Let a time-like vector $P$ be an eigenvector of an operator $p,$
$pP=P,$ and let a time-like vector $Q$ be an eigenvector for $q,$ $qQ=Q.$
Then, $\varpi(p,q)R=(-P\cdot R)\varpi(P,Q),$ is space-like vector,
$\forall\;R.$ The space-like binary relative-velocity, $\varpi(P,Q)=i_{gP}\frac{P\wedge Q}{P\cdot Q},$ was used by Minkowski (1908). Hestenes in 1974 introduced, within a Clifford-geometric algebra, binary relative velocity of $Q$ relative to $P,$ as a Minkowski bivector $\frac{P\wedge Q}{P\cdot Q}.$ In [\'Swierk 1988], we observed that binary relative-velocity $\varpi(P,Q)$ can not parameterize the isometric Lorentz boost. The same definition was introduced independently by many Authors, Matolcsi (1993, \S 4.2-4.3, p. 191), Bini, Carini and Jantzen (1995), Gottlieb (1996), Matolcsi and Goher (2001, p. 89,
Definition (18)), Mitskievich (2005, formula (4.19) on p. 16).

The addition law of the binary space-like relative-velocities were derived by many Authors, (\'Swierk 1988, Matolcsi 1993, \S 4.3; Bini, Carini \& Jantzen 1995). In Oziewicz (2005) we noted that this addition is associative, and compared it with non-associative addition of ternary relative velocities parameterizing the isometric Lorentz boost. Here we show how this associative addition follows from the associativity of the algebra $\obs(\varpi).$

\begin{othm}[Associative addition of binary relative-velocities]\label{add} See
\emph{Figure 3}. Set $\uel\equiv\varpi(p,q),$ $\emph{source}(\uel)=p,$ and $\vel\equiv\varpi(q,r).$ We adhere to the arabic convention of the composition of morphisms, read from the right to the left,
\begin{equation}\varpi(q,r)\circ\varpi(p,q)=\varpi(p,r)=\vel\circ\uel.\end{equation}
The binary composition of composable arrows is associative and has
the following $p$-dependent explicit form,
\begin{equation}(1-\textstyle{\frac{\vel\cdot\uel^{-1}}{c^2}})(\vel\circ\uel)
=\uel+(1-\textstyle{\frac{\uel^2}{c^2}})\vel
p+\frac{1}{c}(\vel\cdot\uel^{-1})p.\label{addit}\end{equation}\end{othm}

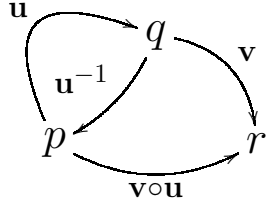
\begin{figure}[h]{\Large$$\xymatrix{&q\ar@/^1pc/[dr]^{\vel}\ar@/^/[dl]_{\uel^{-1}}\\
p\ar@/^3pc/[ur]^{\uel}\ar@/_1pc/[rr]_{\vel\circ\uel}&&r}$$}\caption{Associative
composition of binary velocities-morphisms.\label{com}}\end{figure}

\noindent\textit{Proof.} The velocity-addition follows from the
associativity of an algebra $\obs(\varpi),$ $r(qp)=(rq)p,$ and from
algebra-axiom of Sec. 4, formula (\ref{t1}),
\begin{eqnarray}r(qp)=\trace(qp)\trace(rp)\{\varpi(p,r)+p\}+\trace(qp)r\varpi(p,q),
\nonumber\\(rq)p=\trace(qp)\trace(rq)\{\varpi(p,q)+p\}+\trace(rq)\varpi(q,r)p,\end{eqnarray}
\begin{eqnarray}\trace(qp)\trace(rp)\varpi(p,r)=\trace(qp)\trace(rq)\varpi(p,q)
+\trace(r,q)\varpi(q,r)p\hspace{2cm}\nonumber\\\hspace{5cm}+\trace(qp)
\trace(r(q-p))p-\trace(qp)r\varpi(p,q).\end{eqnarray} The multiplication of an object $r$ with an arrow $\varpi(p,q),$ as given by formula (\ref{t2}), is expressed in terms of $\{\varpi(p,r)+p\}.$ This gives
\begin{eqnarray}\varpi(p,r)=\varpi(q,r)\circ\varpi(p,q)\hspace{6cm}\nonumber\\\hspace{2cm}=
\textstyle{\frac{\trace(pqr)}{\trace(pr)}\left(\varpi(p,q)+
\frac{1}{\trace(pq)}\varpi(q,r)p\right)+c\left(\frac{\trace(pqr)}
{\trace(pr)}-1\right)p.}\label{comp}\end{eqnarray} We need to use
Section 5, expression (\ref{angle}). The associativity of addition
follows from considering a system of four objects/massive bodies,
$p(qrs)=(pqr)s.$ $\Box$\smallskip 

Ungar observed that the addition of (ternary) velocities within isometric special
relativity (within the Lorentz relativity group), is non-associative (Ungar 1988, 2001). In contrast, the addition of binary velocities within relativity groupoid, Theorem \ref{add} and expression (\ref{addit}), is object-dependent (depends on
observer $p$ that is the source of the velocity-arrow $\uel$), however this addition operation is associative.

If $\varpi(p,q)$ is a velocity-operator of $q$ relative to $p,$ then the $\circ$-inverse velocity-morphism, $\varpi(q,p)=(\varpi(p,q))^{-1},$
is a velocity of $p$ relative to $q,$
\begin{eqnarray}\varpi(p,q)\circ\varpi(q,p)=\zel_q\quad\neq\quad\zel_p=
\varpi(p,p)=\varpi(q,p)\circ\varpi(p,q)\nonumber\\\hspace{3cm}=\trace(pq)
\varpi(p,q)+\varpi(q,p)p+c\{\trace(pq)-1\}p.\label{inverse}\end{eqnarray}
Equation (\ref{inverse}) follows from (\ref{comp}), and therefore we have
\begin{equation}\{\varpi(q,p)\}^{-1}=\varpi(p,q)=-\frac{\varpi(q,p)p}{\trace(pq)}-
c\left\{1-\frac{1}{\trace(pq)}\right\}p.\end{equation}

\section{Galilean addition}\setcounter{equation}{22}

Two objects, $p$ and $q,$ possess the same simultaneity iff $pq=p,$ and hence $\trace(pq)=1.$ Therefore in the Galilean algebra of observers the trace of arbitrary string of objects must be $\trace(pq\ldots)=1.$ The Galilean algebra of observers is presented on idempotent-objects, and on the following relations, not independent,
\begin{eqnarray}qp=q=p+\varpi(p,q),\qquad\varpi(p,q)\varpi(r,s)=0,\nonumber\\
q\varpi(p,r)=0,\qquad\varpi(q,r)p=\varpi(q,r).\end{eqnarray} Hence Section 4 with relation (\ref{m}) must be replaced by the reciprocal binary relative-velocity $\varpi(p,q)=q-p,$ \ie relative velocity in Galilean relativity is exactly the \textit{skew} symmetric function of his arguments,
\begin{equation}\varpi(q,r)\circ\varpi(p,q)=\varpi(p,q)+\varpi(q,r)=
(q-p)+(r-q)=\varpi(p,r).\end{equation} Therefore,
$\lim_{c\rightarrow\infty}(\vel\circ\uel)=\vel+\uel,$ because $\vel
p=\vel.$

\section*{Acknowledgments} An author is a member of Sistema Nacional de Investigadores in M\'exico, Expediente \# 15337.

I am grateful to William Page (Kingston, Canada) for insightful observations and valuable discussions,
and to Garret Sobczyk for clarifying and discussing the Hestenes formalism.

\end{document}